\documentclass[12pt]{article}
\usepackage{latexsym,amsmath,amssymb}
\usepackage{amsfonts}
\usepackage{rotating}
\usepackage{multicol}

\usepackage{subfig}

\usepackage[square,numbers,sort&compress]{natbib}
\numberwithin{equation}{section} \numberwithin{table}{section}

\def\half{\mbox{$\frac{1}{2}$}}

\newcommand{\vep}{\mbox{$\varepsilon$}}
\newcommand{\pa}{\mbox{$\partial$}}

\newcommand{\DE}{differential equation}

\newcommand{\PDE}{partial differential equation}

\def\half{\mbox{$\frac{1}{2}$}}

\newcommand{\vslash}{\mbox{\,\rule[-0.37cm]{0.020cm}{1.0cm}\,}}

\newcommand{\Ad}{\mbox{${\rm Ad}\,$}}

\begin{document}

\large

\begin{flushleft}
{\large \bf Group classification and symmetry reductions of a nonlinear Fokker-Planck equation based on the Sharma-Taneja-Mittal entropy
}\\[3 mm]


{ Winter Sinkala\\[3 mm]}
Department of Mathematical Sciences and Computing\\ Faculty of Natural Sciences\\
Walter Sisulu University\\ Private Bag X1, Mthatha 5117, Republic of
South Africa\\ (E-mail addresses: wsinkala@wsu.ac.za).
\end{flushleft}

\vspace{-5mm}

\begin{abstract}
\noindent A nonlinear Fokker-Planck equation that arises in the framework of statistical mechanics based on the Sharma-Taneja-Mittal entropy is studied via Lie symmetry analysis. The equation describes kinetic processes in anomalous mediums where both super-diffusive and subdiffusive mechanisms arise contemporarily and competitively. We perform complete group classification of the equation based on two parameters that characterise the underlying Sharma-Taneja-Mittal entropy. For arbitrary values of the parameters the equation is found to admit a two-dimensional symmetry Lie algebra. We identify and catalogue all the cases in which the equation admits additional symmetries.
We also perform symmetry reductions of the equation and obtain second-order ODEs that describe all essentially different invariant solutions of the Fokker-Planck equation.
\end{abstract}

\smallskip

{\small KEY WORDS: Nonlinear Fokker-Planck equation, Sharma-Taneja-Mittal entropy, Invariant solutions, Lie symmetry analysis, Group classification.}


\section{Introduction}\label{sec1}

Adaptations of the Fokker-Planck equation serve as mathematical models of various problems that arise in physical and biological sciences (see \cite{Caldas14} and the references therein). In statistical mechanics, for example, the Fokker-Planck equation is used to describe kinetic processes in anomalous mediums. In this application the linear version of the equation is considered appropriate for
the description of a wide variety of physical phenomena characterized by short-range interactions and/or short-time
memories, typically associated with normal diffusion. The nonlinear Fokker-Planck equation, on the other hand, is associated with anomalous diffusion, generally associated to non-Gaussian distributions \cite{Scarfone09}.

In this paper we study a nonlinear Fokker-Planck equation (derived in Scarfone and Wada \cite{Scarfone09}) in the framework of statistical mechanics based on a two-parameter entropy known as the Sharma-Taneja-Mittal (STM) entropy. The equation is a $(1+1)$-\PDE, namely
\begin{equation}\label{nlfpe1}
\frac{\partial u}{\partial t} - \frac{\partial}{\partial x}\big(x u\big) -
  \Omega\, \frac{\partial^2}{\partial x^2}\left[u^{1 +r}\left(\frac{r + \kappa}{2\, \kappa}\,u^{\kappa} - \frac{r
- \kappa}{2\, \kappa}\,u^{-\kappa}\right)\right] = 0,
\end{equation}
where $\Omega$ is a constant diffusion coefficient, $r$ and $\kappa$ are \emph{deformation} parameters that characterise the underlying entropy and essentially define the family of equations (\ref{nlfpe1}). The dependent variable $u\equiv u(x,t)$ is the normalized density distribution describing a
conservative particle system in the velocity-time space $(x,t)$. We shall refer to this equation as the Sharma-Taneja-Mittal nonlinear Fokker-Planck equation (STM-NFPE). For a detailed account of the context and derivation of the equation the interested reader may consult Scarfone and Wada \cite{Scarfone09} wherein Equation (\ref{nlfpe1}) is studied, along with two well-known special cases, via Lie symmetry analysis. Admitted Lie point symmetries are found and used to construct invariant solutions.

The invariant solutions of the STM-NFPE reported in \cite{Scarfone09} are based on two infinitesimal symmetry generators
\begin{equation}\label{sw32}
      X_1 = \pa_t, \quad    X_2 = e^{-t}\pa_x,
\end{equation}
the only ones admitted by the whole family of equations represented by (\ref{nlfpe1}) for arbitrary $r$ and $\kappa$.
  In this paper, we determine all instances, depending on specifications of $r$ and $\kappa$, under which Equation (\ref{nlfpe1}) admits additional symmetries besides those in (\ref{sw32}). This is called complete group classification of (\ref{nlfpe1}) (see \cite{olver,OvsiannikovLV82,BlumanKumei} for more on the subject).

Numerous papers have been devoted to group classification of diverse differential equations. In particular, many classes of nonlinear evolution equations depending on arbitrary functions of one, or at most two, variables have been studied via the method of  group classification. For example, Yung \emph{et al} \cite{Yung94}, Ivanova \emph{et al} \cite{IvanovaI, IvanovaII, Popovych} and Vaneeva \emph{et al} \cite{VaneevaII, VaneevaI} report a series of works on group classification of various forms of nonlinear variable coefficient diffusion-convection equations. There are many other papers that have appeared over the years on group classification of differential equations \cite{Muatjetjeja13, Patsiuk15, Bozhkov14, Ellabany02,Zhdanov99, Sinkala08a}.

In the case of the equation under consideration in this paper, (\ref{nlfpe1}), we have determined via the method of group classification that there are basically four cases corresponding to $r = \pm k$ and $r = -1 \pm k$ in which equation (\ref{nlfpe1}) admits additional symmetries.

The paper is organised as follows. In Section~\ref{sec2}, we introduce elements of Lie symmetry analysis of \DE s. Group classification of (\ref{nlfpe1}) is done in Section~\ref{sec3}. In this section we determine all the instances depending on $r$ and $\kappa$ when the principal Lie algebra of (\ref{nlfpe1}) is extended. In Section~\ref{sec4}, we construct adjoint representations of the symmetry Lie algebras of (\ref{nlfpe1}) corresponding to all the instances when the equation admits additional symmetries. Furthermore, we compute corresponding optimal systems of admitted one-dimensional subalgebras and perform symmetry reductions of (\ref{nlfpe1}).
Finally, we give concluding remarks in Section~\ref{sec5}.

\section{Preliminaries}\label{sec2}

Lie symmetry analysis is one of the most powerful methods for finding analytical solutions of differential equations. It has its origins in studies by the Norwegian mathematician Sophus Lie who began to investigate continuous groups of transformations that leave differential equations invariant. Accounts of the subject and its application to differential equations are covered in many books (see, for example, \cite{Cantwell,Hydon2000,HansStephani} for an introduction and \cite{BlumanKumei,OvsiannikovLV82,BlumanAlexei10,olver} for a more detailed exposition).
Central to methods of Lie symmetry analysis is invariance of a differential equation under a continuous group of transformations.
Consider a one-parameter Lie group of point transformations
\begin{eqnarray}\label{wfm1.2}\nonumber
 \tilde{x} &=& x + \vep\,\xi(x,t,u) + O(\vep^2)\\
 \tilde{t} &=& t + \vep\,\tau(x,t,u) + O(\vep^2)\\\nonumber
 \tilde{u} &=& u + \vep\,\eta(x,t,u) +
 O(\vep^2)
\end{eqnarray}
depending on a continuous parameter $\vep$. This transformation is characterised by its its infinitesimal generator,
\begin{equation}\label{wfm1.4}
    X = \xi(x,t,u)\pa_x + \tau(x,t,u)\pa_t+ \eta(x,t,u)\pa_u.
\end{equation}
A general $(1+1)$ \PDE\ with a dependent
variable $u$ and independent variables $(x,t)$,
\begin{equation}\label{2.3}
    \Delta(x,t,u,u_x,u_t,u_{xx},u_{xt},u_{tt}) = 0
\end{equation}
is invariant under (\ref{wfm1.2}) if and only if
\begin{equation}\label{2.4}
    X^{(2)}\Delta = 0\quad \mbox{when}\quad \Delta = 0,
\end{equation}
where $X^{(2)}$ is the second prolongation of $X$ given by
\begin{equation}\label{X2}
 X^{(2)}= X + \eta_i^{(1)} \pa_{u_i} + \eta_{i_1 i_2}^{(2)} \pa_{u_{i_1 i_2}},\quad i_1, i_2 = 1,2,
\end{equation}
with
 \begin{equation}\label{etas}  \eta_i^{(1)}  = D_i \eta - \left(D_i \xi^j\right)u_j, \quad
\eta_{i_1 i_2}^{(2)} = D_{i_2} \eta_{i_1}^{(1)} - \left(D_{i_2}\xi^j\right)u_{i_1 j}, \quad i,i_k,j = 1,2,
\end{equation}
where $u_i=\frac{\partial u}{\partial x^i}$, $u_{i_1 i_2}=\frac{\partial^2 u}{\partial x^{i_1}\partial x^{i_2}}$, $i,i_j=1,2$, $(x^1,x^2) = (x,t)$, $(\xi^1,\xi^2) = (\xi,\tau)$ and $D_i$ denotes the total differential operator with respect to $x^i$:
\begin{equation}\label{Di}
    D_i=\frac{\partial}{\partial x^i}+u_i\frac{\partial}{\partial{u}}+u_{ij}\frac{\partial}{\partial{u_j}} + u_{ijk}\frac{\partial}{\partial u_{jk}}+ \cdots
\end{equation}
The Einstein summation convention is adopted in (\ref{X2}), (\ref{etas}) and (\ref{Di}).
The invariance condition (\ref{2.4}) yields an  over-determined system of linear partial differential
equations (determining equations) for the symmetry group of equation (\ref{2.3}). The infinitesimals $\xi$, $\tau$ and $\eta$ are then determined as a general solution to the determining equations.
 If the infinitesimals contain more than one arbitrary constant, the resulting multi-parameter infinitesimal generator is split into single parameter generators, which constitute a basis for the symmetry Lie algebra of (\ref{2.3}).

When equation (\ref{2.3}) has arbitrary elements (parameters and/or functions) the nature and dimension of the admitted symmetry Lie algebra typically depends on these arbitrary elements. For many equations modelling natural phenomena it is desirable that they possess non-trivial symmetry. This is because in such cases it is possible to obtain a lot of information about solutions of the equation, including reduction of multidimensional equations to ordinary differential equations, constructing classes of exact and approximate solutions. Furthermore, mathematical models must be of such a form that they are consistent with important natural principles of physics such as conservation laws of energy, momentum, etc. It turns out that this beauty of equations of mathematical physics is often encoded in the symmetries admitted by the equation \cite{WIF2002}. From this point of view, therefore, we seek to specify the arbitrary elements in the equation in such a way as to increase the dimension of the admitted symmetry Lie algebra. This is the essence of the group classification method.
Seminal work on group classification was done by Sophus Lie \cite{lie} (see also
\cite{OvsiannikovLV82}) who investigated linear
second-order partial differential equations (PDEs) with
two independent variables.
  For a more detailed account of the method of group classification
of differential equations the reader is refereed to Gazizov and Ibragimov~\cite{Gazizov98} and Olver~\cite{olver} (and the references contained therein).

\section{Lie point symmetries of the STM-NFPE}\label{sec3}
Clearly, when $r = -1/2$ and $k=\pm1/2$  equation (\ref{nlfpe1}) is reduced to a first-order linear equation and therefore admits an infinite dimensional Lie algebra.  These cases will not be considered. For other values of the parameters, suppose
\begin{equation}\label{DeltaOperator}
 X = \xi(x,t,u)\,\pa_x + \tau(x,t,u)\,\pa_t + \eta(x,t,u)\,\pa_u,
\end{equation}
where $\xi$, $\tau$ and $\eta$ are arbitrary functions,
is admitted by Equation (\ref{nlfpe1}). Then the invariance condition dictates the following:
\begin{equation}\label{InvCondBond}
 X^{(2)}\left\{\frac{\partial u}{\partial t} - \frac{\partial}{\partial x}\big(x u\big) -
  \Omega\, \frac{\partial^2}{\partial x^2}\left[u^{1 +r}\left(\frac{r + \kappa}{2\, \kappa}\,u^{\kappa} - \frac{r
- \kappa}{2\, \kappa}\,u^{-\kappa}\right)\right]\right\}\vslash_{(\ref{nlfpe1})} = 0,
\end{equation}
where $X^{(2)}$ is the second-prolongation of $X$ defined in (\ref{X2}).
In the usual way (\ref{InvCondBond}) translates into the following system of  determining equations, thanks to the handiness of Mathematica~\cite{Mathematica}:
\begin{equation}\label{fpdet1}
  {\xi_u} = 0
\end{equation}
\begin{equation}\label{fpdet2}
  {\tau_u} = 0
\end{equation}
\begin{equation}\label{fpdet3}
  {\tau_x} = 0
\end{equation}
\begin{eqnarray}\label{fpdet4}
&& \xi + {\xi_t} - x\,{\xi_x} + x\,{\tau_t} - 2\,\left( k - r \right) \,u^{r - k - 1}\,\varphi(-k,r)\,{\eta_x}\nonumber\\
&&  +\,2\,\left( k + r \right) \,u^{r + k - 1}\,\varphi(k,r)\,{\eta_x} -\,u^{r - k}\,\varphi(-k,r)\,\left( {\xi_{xx}} - 2\,{\eta_{xu}} \right)\nonumber\\
&&  -\, u^{k + r}\,\varphi(k,r)\,\left( {\xi_{xx}} - 2\,{\eta_{xu}} \right) = 0,
\end{eqnarray}
\begin{eqnarray}\label{fpdet5}
\eta + u\,{\tau_t} - {\eta_t} - u\,{\eta_u} + x\,{\eta_x} + u^{r - k}\,\varphi(-k,r)\,{\eta_{xx}} +
  u^{k + r}\,\varphi(k,r)\,{\eta_{xx}} = 0,
\end{eqnarray}
\begin{eqnarray}\label{fpdet6}
\eta\,\left( r + k - 1 \right) \,\left( k + r \right) \,u^{r + k - 1}\,\varphi(k,r) + \eta\,u^{r - k - 1}\,\vartheta(k,r) \nonumber\\
   +\,
  \left[ \left( r - k \right) \,u^{r - k}\,\varphi(-k,r) + \left( k + r \right) \,u^{k + r}\,\varphi(k,r) \right] \,
   \left({\tau_t} + {\eta_u} -2\,{\xi_x}\right)  \nonumber\\
   +\, \left( u^{1 + k + r}\,\varphi(k,r) + \frac{\Omega\,u^{1 - k + r}\,\vartheta(k,r)}{2\,k\,\varphi(k,-r)}
     \right) \,{\eta_{{uu}}} = 0,
\end{eqnarray}
\begin{eqnarray}\label{fpdet7}
\eta\,\left( r - k \right) \,u^{r - k - 1}\,\varphi(-k,r) + \eta\,\left( k + r \right) \,u^{r + k - 1}\,\varphi(k,r)~~~~~~~~ \nonumber\\
   -\,u^{r - k}\,\varphi(-k,r)\,\left( 2\,{\xi_x} - {\tau_t} \right) - u^{k + r}\,\varphi(k,r)\,\left( 2\,{\xi_x} - {\tau_t} \right) = 0,
\end{eqnarray}
where
\[ \varphi(k,r) = \frac{\Omega\,\left( k + r \right) \,\left( 1 + k + r \right) }{2\,k}, \quad \vartheta(k,r) = \frac{\Omega\,\left[1 - {\left( k - r \right) }^2 \right] \,{\left( k - r \right) }^2}{2\,k}.\]
From the equations (\ref{fpdet1})--(\ref{fpdet3}) and (\ref{fpdet7}), we obtain that
\begin{eqnarray}
  \xi(x,t,u) &=& \xi(x,t) \label{f1} \\
  \tau(x,t,u) &=& \tau(t) \label{f2}\\
  \eta(x,t,u) &=& \frac{u \left[ \varphi(-k, r) + \varphi(k, r)\,u^{2\,k} \right]
    \left(2\,{f}_x  - {g}_t\right) }{(r - k)\varphi(-k, r) + (k + r)\varphi(k, r)\,u^{2\,k}}.\label{f3}
\end{eqnarray}
The outstanding equations, (\ref{fpdet4})--(\ref{fpdet6}), are now expressed in terms of and solved for $\xi$ and $\tau$. Equation (\ref{fpdet6}) in particular becomes
\begin{equation}\label{3.13}
  \frac{\left[u^{k + r}\,\psi_1(k,r) + u^{3\,k + r}\,\psi_2(k,r) +
  u^{5\,k + r}\,\psi_3(k,r)\right]\left[\tau^\prime - 2\,\xi_x\right]}{\psi_0(k,r,u)} = 0,
\end{equation}
where $\prime$ denotes the differentiation with respect to $t$, and
\begin{eqnarray*}
  \psi_0 &=& {\left[ \left(k - r - 1 \right) \,{\left( k - r \right) }^2 +
     {\left( k + r \right) }^2\,\left( 1 + k + r \right) \,u^{2\,k} \right] }^3 \\
  \psi_1 &=& 2\,k\,\left( 1 + 2\,k \right) \,\Omega\,{\left( k - r \right) }^3\,{\left( 1 - k + r \right) }^2
  \\
  && \times\,\left[ k^4 + r^2\,{\left( 1 + r \right) }^2 - k^2\,\left( 1 + 2\,r + 2\,r^2 \right)  \right] \\
  \psi_2 &=& 4\,k\,\Omega\,\left[k^4 + r^2\,{\left( 1 + r \right) }^2 - k^2\,\left( 1 + 2\,r + 2\,r^2 \right)  \right] \\  && \times\,
  \Big[r^3\,{\left( 1 + r \right) }^2 -2\,k^6 + k^4\,\left( 2 + 5\,r + 4\,r^2 \right)\\
&&~~~~~~~~~~ - \, k^2\,r\,\left( 1 + 4\,r + 6\,r^2 + 2\,r^3 \right)  \Big]\\
 \psi_3  &=& 2\,k\,\left( 2\,k - 1\right) \,\Omega\,{\left( k + r \right) }^3\,{\left( 1 + k + r \right) }^2
  \\
  && \times\,\left[ k^4 + r^2\,{\left( 1 + r \right) }^2 - k^2\,\left( 1 + 2\,r + 2\,r^2 \right)  \right].
\end{eqnarray*}
Clearly equation (\ref{3.13}) is solved if and only if $\psi_0(k,r,u)\ne 0$ and
\begin{equation}\label{3.14}
    \tau^\prime - 2\,\xi_x = 0
\end{equation}
or
\begin{equation}\label{3.15}
     u^{k + r}\,\psi_1(k,r) + u^{3\,k + r}\,\psi_2(k,r) +
  u^{5\,k + r}\,\psi_3(k,r) = 0.
\end{equation}
In the light of (\ref{f1})--(\ref{f3}), the solution of (\ref{3.14}) together with (\ref{fpdet4}) and (\ref{fpdet5}) leads to
\begin{equation}\label{}
    \xi= \vep_1\,e^{-t}, \quad \tau = \vep_2, \quad \eta = 0,
\end{equation}
from which we obtain the infinitesimal symmetry generators $X_1$ and $X_2$ in (\ref{sw32}). This means that the principal Lie algebra of (\ref{nlfpe1})  is spanned by $X_1$ and $X_2$.
Instances in which (\ref{nlfpe1}) admits additional symmetries are obtained from solutions of (\ref{3.15}). Solving this equation for $r$, by setting $\psi_1(k,r) = \psi_2(k,r) = \psi_3(k,r) = 0$, we obtain
\begin{equation}\label{cases}
    r = \pm k \quad \mbox{or}  \quad  r = -1 \pm k.
\end{equation}
We consider each of these cases in turn to determine symmetries admitted by (\ref{nlfpe1}) in these instances. This essentially means solving the outstanding determining equations (\ref{fpdet4}) and (\ref{fpdet5}) for $\xi$ and $\tau$. For each of the parameter specifications in (\ref{cases}) the infinitesimal symmetry generators admitted by (\ref{nlfpe1}) include $X_1$ and $X_2$ specified in (\ref{sw32}). Additional ones are admitted in the various cases as presented below in Table~\ref{TB1}.

\renewcommand\arraystretch{1.50}
\begin{table}[!h]
\caption{Infinitesimal symmetry generators of Eqn~(\ref{nlfpe1}) in the Cases $r=\pm k$ and $r=-1\pm k$.}\label{TB1}
\vspace{-6mm}
\[\begin{array}{c|c|l}
\hline
  \mbox{Case} & \mbox{Specifications of $r$ and $k$~~~~~~~~~~~} & \mbox{Admitted infinitesimal generators}\\\hline\hline
\mbox{A}  & \begin{array}{cl}
                \mbox{(i)} & r=k,~ k = -\frac{2}{3} \\
                \mbox{(ii)} & r=-k,~ k = \frac{2}{3} \\
                \mbox{(iii)} & r=-1+k,~k = -\frac{1}{6} \\
                \mbox{(iv)} &r= -1-k,~ k = \frac{1}{6}~~~~~~~~~~~
              \end{array}
  & \begin{array}{l}
X_1 = \pa_t, \quad    X_2 = e^{-t}\pa_x\\
   X_3 = e^{-2\,t/3} \left(x\,\pa_x - \pa_t - u\,\pa_u\right)\\ X_4 = x\,e^t\left(x\,\pa_x - 3\,u\,\pa_u\right)\\X_5 = x\,\pa_x - \pa_t - \frac{3}{2}\,u\,\pa_u
                                \end{array}\\\hline
\mbox{B}  & \begin{array}{cl}
                \mbox{(i)} & r=k,~ k = -1 \\
                \mbox{(ii)} & r=-k,~ k = 1 \\
                \mbox{(iii)} & r=-1+k,~k = -\frac{1}{2} \\
                \mbox{(iv)} &r= -1-k,~ k = \frac{1}{2}~~~~~~~~~~~
              \end{array} &
                                \begin{array}{l}
X_1 = \pa_t, \quad    X_2 = e^{-t}\pa_x\\
X_3 = x\,\pa_x-\pa_t - u\,\pa_u,\\ X_4 = t\,x\,\pa_x-t\,\pa_t - u\left(t+\mbox{$\frac{1}{2}$}\right)\,\pa_u
                                \end{array}\\\hline
\mbox{C}  & \begin{array}{cl}
                \mbox{(i)} & r=k,~ k\notin \left\{-\half,-\frac{2}{3},-1\right\}\\
                \mbox{(ii)} & r=-k,~ k\notin \left\{\half,\frac{2}{3},1\right\} \\
                \mbox{(iii)} & r=-1+k,~k\notin \left\{\pm\half,-\frac{1}{6}\right\} \\
                \mbox{(iv)} &r= -1-k,~k\notin \left\{\pm\half,\frac{1}{6}\right\}
              \end{array}
 &
                                \begin{array}{l}
X_1 = \pa_t, \quad    X_2 = e^{-t}\pa_x\\
X_3 = e^{-2\,\left( 1 + \delta \right) \,t}\left( x\, \pa_x -\pa_t - u\, \pa_u  \right)\\
X_4 = x\,\pa_x-\pa_t + (u/\delta)\,\pa_u,\\
\mbox{where}~\delta = \left\{
             \begin{array}{cl}
               k & \hbox{for (i)} \\
               -k & \hbox{for (ii)} \\
               k-\half & \hbox{for (iii)} \\
               -k-\half & \hbox{for (iv).}
             \end{array}
           \right.
                                \end{array}
                              \\\hline
\end{array}\]

\end{table}

The commutator tables for the corresponding Lie algebras are given in Tables~\ref{CTa}, \ref{CTb} and \ref{CTc}.

\begin{table}[!h]
\caption{Lie Brackets for Case A.}\label{CTa}
\centering
    \begin{tabular}{c|ccccc}
      \hline
      $[X_i,X_j]$ & $X_1$ & $X_2$ & $X_3$ & $X_4$ & $X_5$  \\\hline\hline
      $X_1$ & $0$ & $-X_2$ & $-\frac{2}{3}X_3$ & $X_4$ & $0$  \\
      $X_2$ & $X_2$ & $0$ & $0$ & $2\left(X_1 + X_5\right)$ & $0$  \\
      $X_3$ & $\frac{2}{3}X_3$ & $0$ & $0$ & $0$ & $-\frac{2}{3}X_3$  \\
      $X_4$ & $-X_4$ & $-2\left(X_1 + X_5\right)$ & $0$ & $0$ & $0$  \\
      $X_5$ & $0$ & $0$ & $\frac{2}{3}X_3$ & $0$ & $0$  \\
      \hline
    \end{tabular}
\end{table}
\begin{table}[!h]
\caption{Lie Brackets for Case B.}\label{CTb}
\centering
    \begin{tabular}{c|cccc}
      \hline
      $[X_i,X_j]$ & $X_1$ & $X_2$ & $X_3$ & $X_4$ \\\hline\hline
      $X_1$ & $0$ & $-X_2$ & $0$ & $X_3$ \\
      $X_2$ & $X_2$ & $0$ & $0$ & $0$ \\
      $X_3$ & $0$ & $0$ & $0$ & $-X_3$ \\
      $X_4$ & $-X_3$ & $0$ & $X_3$ & $0$ \\
      \hline
    \end{tabular}
\end{table}

\begin{table}[!h]
\caption{Lie Brackets for Case C, where $\delta = k, -k, k-\half, -k-\half$ in the specifications (i), (ii), (iii) and (iv), respectively.}\label{CTc}
\centering
    \begin{tabular}{c|cccc}
      \hline
      $[X_i,X_j]$ & $X_1$ & $X_2$ & $X_3$ & $X_4$ \\\hline\hline
      $X_1$ & $0$ & $-X_2$ & $-2 (1 + \delta)X_3$ & $0$ \\
      $X_2$ & $X_2$ & $0$ & $0$ & $0$ \\
      $X_3$ & $2 (1 + \delta)X_3$ & $0$ & $0$ & $-2 (1 + \delta)X_3$ \\
      $X_4$ & $0$ & $0$ & $2 (1 + \delta)X_3$ & $0$ \\
      \hline
    \end{tabular}
\end{table}

\section{Invariant solutions of the STM-NFPE}\label{sec4}

Each of the infinitesimal symmetry generators admitted by (\ref{nlfpe1}) can be used to construct a family of invariant solutions of the equation.
A function $u = \Theta(x, t)$ is an invariant solution of (\ref{2.3}) arising from $X$ if it is a solution of (\ref{2.3}) and satisfies the \emph{invariant surface} condition,
\begin{equation}\label{4.1}
  X\left(u - \Theta(x, t)\right)=0\quad\mbox{when}\quad u = \Theta(x, t).
\end{equation}
The construction of invariant solutions proceeds in a very algorithmic fashion. For each infinitesimal symmetry generators $X$, one determines from solutions of the associated corresponding system,
\begin{equation}\label{}
    \frac{d x}{\xi} = \frac{d t}{\tau} = \frac{d u}{\eta},
\end{equation} two independent invariants $r(x, t, u)$ and $v(x, t, u)$ (with $v_u\ne 0$) of the associated group.
The form of the invariant solution arising from $X$ is now obtained from $v=F(r)$ or
\begin{equation}\label{4.3}
   u=\Theta(x,t)
\end{equation}
when solved for $u$. Upon substitution of (\ref{4.3}) into (\ref{2.3}) we obtain an ODE that defines $\Theta$, the solution of which completes the construction of the invariant solution.

To avoid ``duplicating" invariant solutions we determine optimal systems (in the usual way \cite{olver}). For each of the Lie algebras represented in Table~\ref{TB1}, we construct an adjoint representation of the underlying Lie group via the Lie series
\begin{equation}\label{Lie-series}
\mbox{Ad}\big(\exp(\vep\,X_i)\big)X_j = X_j - \vep[X_i,X_j] +
\frac{1}{2}\,\vep^2[X_i,[X_i,X_j]] - \cdots,
\end{equation}
where $[X_i,X_j]$ is the Lie bracket of $X_i$ and $X_j$. For the cases identified in Table~\ref{TB1} the adjoint representations are presented in Table~\ref{adta}, Table~\ref{adtb} and Table~\ref{adtc} for Cases A, B and C, respectively, where the $(i,j)$-th entry indicates
 $\Ad(\exp(\vep X_i))X_j$.

\begin{table}[!h]
\caption{Adjoint representations for Case A.}\label{adta}
\vspace{-5mm}
\[ \begin{array}{c|ccccc}
      \hline
      \mbox{Ad} & X_1 & X_2 & X_3 & X_4 & X_5  \\\hline\hline
 X_1 & X_1 & e^{\vep }X_2 & e^{\frac{2}{3}\,\vep }X_3 & e^{-\vep }X_4 & X_5 \\
 X_2 & X_1 - \vep\,X_2  & X_2 & X_3 & {\vep }^2 X_2  +  X_4  - 2\,\vep(X_1 + X_5) & X_5 \\
   X_3 & X_1 - \frac{2\,\vep }{3}\,X_3 & X_2 & X_3 & X_4 & X_5 + \frac{2\,\vep }{3}\,X_3  \\
   X_4 & X_1 + \vep\,X_4  &  X_2 + {\vep }^2 X_4 +
   2\,\vep(X_1 + X_5) & X_3 & X_4 & X_5 \\
 X_5 & X_1 & X_2 & e^{-\frac{2 }{3}\,\vep}X_3 & X_4 & X_5 \\
      \hline
    \end{array}\]
\end{table}

\begin{table}[!h]
\caption{Adjoint representations for Case B.}\label{adtb}
\[
    \begin{array}{c|cccc}
      \hline
      \mbox{Ad}& X_1 & X_2 & X_3 & X_4\\\hline\hline
 X_1 &  X_1 & e^{\vep }X_2 & X_3 & X_4 - \vep\,X_3  \\
 X_2 & X_1 - \vep\,X_2  & X_2 & X_3 & X_4 \\
 X_3 & X_1 & X_2 & X_3 & X_4 + \vep\,X_3  \\
 X_4 & X_1 +  \left( 1 - e^{-\vep } \right)X_3\  & X_2 & e^{-\vep }\,X_3 & X_4 \\
      \hline
    \end{array}\]
\end{table}

\begin{table}[!h]
\caption{Adjoint representations for Case C.}\label{adtc}
\[
    \begin{array}{c|cccc}
      \hline
    \mbox{Ad} & X_1 & X_2 & X_3 & X_4 \\\hline\hline
X_1 & X_1 & e^{\vep }X_2 & e^{2\,\left( 1 + \delta \right) \,\vep }X_3 & X_4 \\
X_2 & X_1 - \vep\,X_2\  & X_2 & X_3 & X_4 \\
X_3 & X_1 - 2\,\left( 1 + \delta \right) \,\vep\,X_3  & X_2 & X_3 & X_4 + 2\,\left( 1 + \delta \right) \,\vep\,X_3  \\
X_4 & X_1 & X_2 & e^{-2\,\left( 1 + \delta \right) \,\vep }X_3 & X_4 \\
      \hline
    \end{array}\]
\end{table}
Following Olver's approach \cite{olver}, we use the adjoint representations to identify equivalent infinitesimal symmetry generators. The results are presented in Table~\ref{NFPEOS}.

\begin{table}[h]
\caption{Optimal system of one-dimensional subalgebras of the STM-NFPE; $\alpha$, $\beta$ and $\gamma$ are arbitrary constants.}\label{NFPEOS}
\[
    \begin{array}{c|c}
      \hline
\mbox{Case A} & \begin{array}{cl}
                  \mbox{(a)} & \beta X_2 + X_4 + \gamma X_5 \\
                  \mbox{(b)} & \alpha X_1 + \beta X_2 + X_4 \\
                  \mbox{(c)} &  X_1 + X_3 \\
                  \mbox{(d)} & X_1
                \end{array}
  \\
      \hline
    \end{array}~~~~~~~~\begin{array}{c|c}
      \hline
\begin{array}{c}
  \mbox{Case B} \\[-1ex]
  \&\\[-1ex]
  \mbox{Case C} \\
 ~
\end{array}
 & \begin{array}{cl}
                  \mbox{(a)} & \alpha X_1 + X_4 \\
                  \mbox{(b)} & \alpha X_1 + X_3 \\
                  \mbox{(c)} &  X_1\\
                             &
                \end{array}
  \\
      \hline
    \end{array}
\]
\end{table}

Invariant solutions of (\ref{nlfpe1}) arising from the representative infinitesimal symmetry generators in the optimal systems represented in Table~\ref{NFPEOS} are obtained in the usual way.
\begin{enumerate}
 \item[\textbf{1.}] Case A.
  \begin{itemize}
    \item $\beta X_2 + X_4 + \gamma X_5$
\begin{equation}\label{}
    u(x,t)=e^{\frac{3\,t}{2}} {\left( 1 + x^2 e^{2\,t}/{\beta} \right) }^{-\frac{3}{2}} y\left(\zeta\right),
\end{equation}
\[\zeta = {\sqrt{\beta}}\,t/\gamma + \tan^{-1}\left(e^t\,x/\sqrt{\beta}\right), \quad \beta\, \gamma\ne 0,\]
where
\begin{equation}\label{}
18\,\beta^{\frac{3}{2}}\,y^{\frac{7}{3}}\,y^\prime + 9\,\beta\,\gamma\,y^{\frac{10}{3}} + 2\,\theta\,\gamma\,\Omega\,\left[3\,y\,y^{\prime\prime} - 9\,y^2 - 4\,(y^\prime)^2 \right] = 0.
\end{equation}
with
\begin{equation}\label{EqTheta1}
 \left\{
    \begin{array}{ll}
      \theta = 1, & \hbox{for Parts (i) and (ii),} \\
      \theta = 4, & \hbox{for Parts (iii) and (iv).}
    \end{array}
  \right.
\end{equation}
    \item $\alpha X_1 + \beta X_2 + X_4$
\begin{equation}\label{}
    \!\!\!\!\!\!\!\!\!\!\!\!\!\!  u(x,t)=\frac{8\,\alpha^3\,{\omega }^3\,e^{\frac{3 \left( 1 - \omega  \right)t }{2}} y\left(\zeta\right)}{{\left[ 2\,e^t\,x + \alpha(1+\omega)  \right] }^3},
 \quad \zeta = \frac{\alpha\,\left( \omega - 1  \right) - 2\,e^t\,x }{\omega \,\left[ 2\,e^t\,x + \alpha(1+\omega)  \right]e^{\omega\,t } },
\end{equation}
where
 \begin{equation}\label{}
  18\,{\omega }^5\,\alpha^2\,\zeta\,y^{\frac{7}{3}}\,y^{\prime} + \rho\,y^{\frac{10}{3}} + 2\,\theta\,\Omega\left[4\,(y^{\prime})^2 - 3\,y\,y^{\prime\prime}\right] = 0,
 \end{equation}
with $\theta$ as defined in (\ref{EqTheta1}), and
\[\omega = \sqrt{1 - 4\,\beta/\alpha^2}, \quad \alpha^2 - 4\,\beta \ge 0\quad \mbox{and}\quad \rho = 9\,\alpha^2 {\omega }^4\,\left(3\,\omega - 1 \right).\]
    \item $X_1 + X_3$
\begin{equation}\label{}
    u(x,t)=\frac{e^t\,y(\zeta)}{{\left( e^{\frac{2\,t}{3}} - 1 \right) }^{\frac{3}{2}}}, \quad \zeta = \frac{e^t\,x}{{\left( e^{\frac{2\,t}{3}} - 1 \right) }^{\frac{3}{2}}},
\end{equation}
where
\begin{equation}
    \theta\left[4\,\Omega\,(y^{\prime})^2 - 3\,\Omega\,y\,y^{\prime\prime}\right] + 9\,\zeta\,y^{\frac{7}{3}}\,y^{\prime} + 9\,y^{\frac{10}{3}} = 0 \label{4.12}
\end{equation}
with $\theta$ as defined in (\ref{EqTheta1}).
    \item $X_1$
\begin{equation}\label{}
    u(x,t)=y(\zeta), \quad \zeta = x,
\end{equation}
where $y$ is a solution of (\ref{4.12}) with $\theta=1$.
  \end{itemize}
 \item[\textbf{2.}] Case B.
  \begin{itemize}
    \item $\alpha X_1 + X_4$
\begin{equation}\label{}
    u(x,t)=e^t\,{\left( t - \alpha \right) }^{\alpha + \half}\,y(\zeta),\quad \zeta = e^t\,{\left( t - \alpha \right) }^{\alpha}\,x,
\end{equation}
where
\begin{equation}\label{}
 2\,\theta\,\Omega\left[y\,y^{\prime\prime} - 2\,(y^{\prime})^2\right] + 2\,\alpha\,\zeta\,y^3\,y^{\prime} + \left( 1 + 2\,\alpha \right) \,y^4 = 0.
\end{equation}
with
\begin{equation}\label{EqTheta}
 \left\{
    \begin{array}{ll}
      \theta = 1, & \hbox{for Parts (i) and (ii),} \\
      \theta = 2, & \hbox{for Parts (iii) and (iv).}
    \end{array}
  \right.
\end{equation}
 \item $\alpha X_1 + X_3$
\begin{equation}\label{}
    u(x,t)=e^{\frac{t}{1 - \alpha}}\,y(\zeta), \quad \zeta = e^{\frac{t}{1 - \alpha}}\,x,
\end{equation}
where
\begin{equation}\label{}
 \theta\,\Omega\,(\alpha - 1)\,[ 2\,(y^{\prime})^2  - y\,y^{\prime\prime}]  + \alpha(y^4 + \zeta \, y^3\,y^{\prime}) = 0.
\end{equation}
with $\theta$ as defined in (\ref{EqTheta}).
    \item $X_1$
\begin{equation}\label{}
    u(x,t) = y(x),
\end{equation}
where
\begin{equation}\label{}
  \theta\,\Omega\left[y\,y^{\prime\prime} - 2\,(y^{\prime})^2\right] - x\, y^3\,y^{\prime} -y^4 = 0,
\end{equation}
with $\theta$ as defined in (\ref{EqTheta}).
  \end{itemize}

 \item[\textbf{3.}] Case C: (i) \& (ii) \ [$\delta = k$ for (i),\, $\delta =-k$ for (ii)]
  \begin{itemize}
    \item $\alpha X_1 + X_4$
\begin{equation}\label{}
    u(x,t)=e^{\frac{t}{\delta(\alpha -  1)}}\,y(\zeta), \quad \zeta = x \,e^{\frac{t}{1 - \alpha}},
\end{equation}
where
\begin{equation}\label{}
  y\,\left[ \lambda\,y  - \delta\,\alpha\,\zeta\,y^{\prime}  \right]  -
 \Omega\,\delta\,\omega\, y^{2\,\delta } \,\left[ y\,y^{\prime\prime} + 2\,\delta\,(y^{\prime})^2  \right] = 0,
\end{equation}
with
 \[\omega = \left( \alpha - 1\right) \,\left( 1 + 2\,\delta  \right), \quad \lambda=1 + \delta(1  - \alpha).\]
    \item $X_1 + X_3$
\begin{equation}\label{}
 u(x,t)=\frac{e^t\,y(\zeta)}{{\left( e^{2\,\left( 1 + \delta  \right)\,t } - 1 \right) }^{\frac{1}{2(1 +\delta) }}}, \quad \zeta = \frac{e^t\,x}{{\left( e^{2\,\left( 1 + \delta  \right)\,t } - 1 \right) }^{\frac{1}{2(1 +\delta) }}}.
\end{equation}
where
\begin{equation}\label{}
  y^2 + \Omega\,\left( 1 + 2\,\delta  \right) y^{2\,\delta } \left[ y\,y^{\prime\prime} + 2\,(y^{\prime})^2\,\delta  \right]  + \zeta\,y\,y^{\prime}= 0.
\end{equation}
    \item $X_1$
\begin{equation}\label{}
    u(x,t) = y(x),
\end{equation}
where
\begin{equation}\label{}
  \Omega\,\left( 1 + 2\,\delta  \right)\,y^{2\,\delta - 1} \left[ y\,y^{\prime\prime} + 2\,\delta \,(y^{\prime})^2 \right]  + x\,y^{\prime} + y = 0.
\end{equation}
  \end{itemize}
  \item[\textbf{4.}] Case C: (iii) \& (iv)  [$\delta = k$ for (iii) \& $\delta =-k$ for (iv)]
  \begin{itemize}
    \item $\alpha X_1 + X_4$
\begin{equation}\label{}
    u(x,t)=\exp\left\{\frac{2\,t}{\nu\,\left(\alpha -1 \right)}\right\}\,y(\zeta), \quad \zeta = x\,e^{\frac{t}{1 - \alpha}},
\end{equation}
where
\begin{equation}\label{}
  \Omega\,\nu\,y^{1 + 2\,\delta }\,y^{\prime\prime}  + \Omega\,{\nu }^2\,y^{2\,\delta }\,(y^{\prime})^2 + \mu \,\zeta\, y^2\,y^{\prime}  + \omega \,y^3= 0,
\end{equation}
with
\[\mu = \frac{\alpha}{\alpha - 1}, \quad \omega = 1 + \frac{2}{\left( 1 - \alpha \right) \,\nu }, \quad  \nu = 2\,\delta - 1.\]
    \item $X_1 + X_3$
\begin{equation}\label{}
    u(x,t)=\frac{e^t\,y(\zeta)}{{\left[e^{(1 + 2\,\delta)t } - 1\right] }^{\frac{1}{1 + 2\,\delta }}},\quad \zeta = \frac{x\,e^t}{{\left[e^{(1 + 2\,\delta)t } - 1\right] }^{\frac{1}{1 + 2\,\delta }}},
\end{equation}
where
\begin{equation}\label{ws11}
   y^3 + \Omega\,{\rho}^2\,y^{2\,\delta }\,(y^{\prime})^2 - \Omega\,\rho\,y^{1 + 2\,\delta }\,y^{\prime\prime} +
  \zeta\,y^2\,y^{\prime} = 0, \quad \rho = 1 - 2\,\delta.
\end{equation}
    \item $X_1$
\begin{equation}\label{}
    u(x,t) = y(\zeta),\quad \zeta = x,
\end{equation}
where $y$ is any solution of (\ref{ws11}).

  \end{itemize}

\end{enumerate}

\section{Concluding remarks}\label{sec5}
In the work reported in this paper we have carried out complete group classification of the STM-NFPE and determined that the principal Lie algebra of the equation is a two-dimensional Lie algebra. We have further established that the principal Lie algebra extends only if $r=\pm k$ or $r=-1\pm k$ and that in each of these cases the admitted symmetry Lie algebra extends further for specific values of the parameters. We remark here that the identified parameter specifications, for which the principal Lie algebra extends, define particular interesting entropies associated with the STM-NFPE and require further study. In each of these instances the equation is endowed with remarkable physical properties and is amenable to solution via routines of Lie symmetry analysis. As such we have performed symmetry reductions of the STM-NFPE in all such cases, limiting our calculations to essentially different infinitesimal symmetry generators, i.e. those not connected by means of a Lie point transformation of the equation. To do this we have constructed adjoint representations of the symmetry Lie algebras of the STM-NFPE corresponding to all the instances when the equation admits a nontrivial symmetry Lie algebra. Using the adjoint representations, we have determined one-dimensional optimal systems of Lie algebras of the STM-NFPE and used them to perform symmetry reductions of the equation. We have thereby characterised practically all invariant solutions of STM-NFPE by second-order ODEs. Our results have therefore revealed new possibilities for analytical and numerical
studies of the equation.

\section*{Acknowledgements}
The author would like to thank the Directorate of Research Development of Walter
Sisulu University for continued financial support.

%

\end{document}